\documentclass{naco}
\usepackage{amsmath}
  \usepackage{paralist}
  \usepackage{graphics} 
  \usepackage{epsfig} 
\usepackage{graphicx}  \usepackage{epstopdf}
 \usepackage[colorlinks=true]{hyperref}
\hypersetup{urlcolor=blue, citecolor=red}

\def\currentvolume{8}
 \def\currentissue{1}
  \def\currentyear{2018}
   \def\currentmonth{March}
    \def\ppages{1--16}
     \def\DOI{10.3934/naco.2018001}

\newtheorem{theorem}{Theorem}[section]

\theoremstyle{definition}
\newtheorem{definition}[theorem]{Definition}

\title[Linear Reductions of Karp's Problems] 
      {Linearly-growing reductions of Karp's 21 NP-complete problems}

\author[Jerzy A. Filar, Michael Haythorpe and Richard Taylor]{}

\subjclass{Primary: 68Q15; Secondary: 90C10.}
 \keywords{NP-complete, reduction, linear, Karp, complexity, integer programming}

 \email{j.filar@uq.edu.au}
 \email{michael.haythorpe@flinders.edu.au}
 \email{richard.taylor@dsto.defence.gov.au}

\thanks{$^*$ Corresponding author: Michael Haythorpe: michael.haythorpe@flinders.edu.au}

\thanks{The reviewing process was handled by A. Baghirov.}

\begin{document}
\maketitle

\centerline{\scshape Jerzy A. Filar}
\medskip
{\footnotesize
 \centerline{School of Mathematics and Physics}
   \centerline{University of Queensland}
   \centerline{ St Lucia, QLD 4072, Australia}
} 

\medskip

\centerline{\scshape Michael Haythorpe$^*$}
\medskip
{\footnotesize
 \centerline{ College of Science and Engineering}
 \centerline{Flinders University}
   \centerline{Bedford Park, SA 5042, Australia}
}

\medskip

\centerline{\scshape Richard Taylor}
\medskip
{\footnotesize
 \centerline{ Defence Science and Technology Group}
   \centerline{Canberra, ACT 2600, Australia}
}

\bigskip


\begin{abstract}
We address the question of whether it may be worthwhile to convert certain, now classical, NP-complete problems to one of a smaller number of kernel NP-complete problems. In particular, we show that Karp's classical set of 21 NP-complete problems contains a kernel subset of six problems with the property that each problem in the larger set can be converted to one of these six problems with only linear growth in problem size. This finding has potential applications in optimisation theory because the kernel subset includes 0-1 integer programming, job sequencing and undirected Hamiltonian cycle problems.
\end{abstract}

\section{Introduction}

The set of {\em decision problems} is the set of all problems for which any instance has an answer of YES or NO. Among this set is the subset of {\em nondeterministic polynomial time problems}, denoted $\mathcal{NP}$ which are decision problems with the following property. Consider any instance of an $\mathcal{NP}$ problem. If the instance has the answer YES, it should be possible to provide evidence of that answer which can be checked by a deterministic Turing machine in time bounded above by a polynomial function in the input size of the instance. Such evidence is called a {\em certificate}.

One of the most important $\mathcal{NP}$ problems is {\em boolean satisfiability}, which asks whether a set of boolean variables can be assigned values of TRUE or FALSE in order to make a given logical expression evaluate to TRUE. Certainly the problem is in $\mathcal{NP}$, because if the answer is YES, simply providing such an assignment of values for each boolean variable suffices as a certificate. Any logical expression can be written in {\em conjunctive normal form}, that is, it is a conjunction of clauses which, in turn, only involve OR and NOT connectives.

In 1971, the famous Cook-Levin theorem \cite{cook} demonstrated that any $\mathcal{NP}$ problem can be {\em reduced} to boolean satisfiability in conjunctive normal form (SAT) in the following sense. Consider any problem $P \in \mathcal{NP}$. Then there exists a polynomial-time algorithm which accepts, as input, any instance of $P$, say $I_P$, and outputs a new instance $I_B$ of SAT with input size bounded above by a polynomial function in the input size of $I_P$. Then $I_P$ and $I_B$ have the same answer, and if the answer is YES, there exists another polynomial-time algorithm which accepts, as input, any valid certificate of $I_B$, and outputs a valid certificate of $I_P$.

A corollary of the Cook-Levin theorem is that SAT is, in the worst case, at least as difficult, up to polynomial equivalence as any other NP problem. There is another set of problems called {\em NP-hard problems} which is the set of problems that are at least as difficult as the most difficult $\mathcal{NP}$ problem(s). The intersection of $\mathcal{NP}$ and the set of NP-hard problems is called the set of {\em NP-complete problems}, which therefore includes boolean satisfiability, as well as many other problems. If a reduction can be constructed from any problem known to be NP-complete to another problem in $\mathcal{NP}$, then that second problem is hence proved to be NP-complete as well.

\begin{center}\begin{figure}[h!]\includegraphics[scale=0.5]{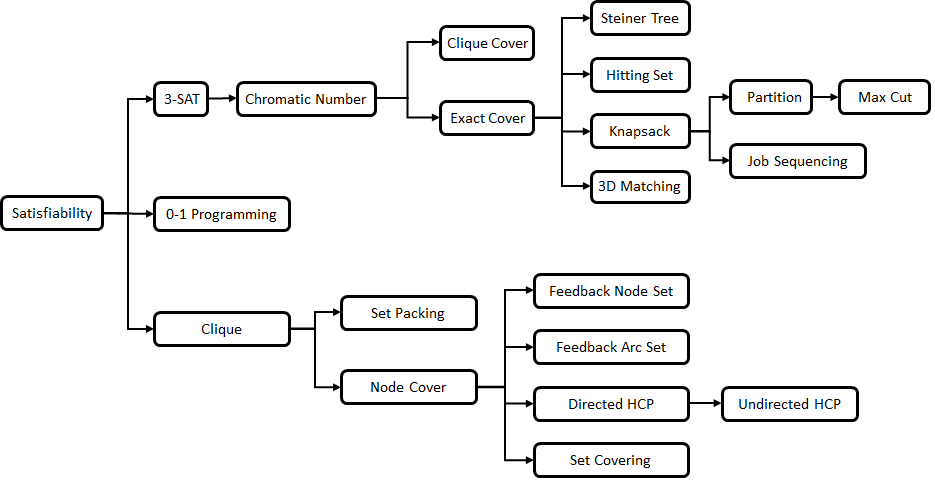}\caption{A tree showing the 21 NP-complete problems identified by Karp, where the edges correspond to individual reductions.}\label{fig-karp}\end{figure}\end{center}

The first major study into this field was by Karp \cite{karp}, who in 1972 provided a list of 21 NP-complete problems (including SAT) by describing twenty such reductions, starting by reducing SAT to three other problems, then reducing those problems to yet more problems, and so forth. A visualisation of the reduction tree is shown in Figure \ref{fig-karp}. We note that Karp's set of 21 problems includes many that have been studied intensively in the context of optimisation, even though in the present context they are cast as decision or feasibility problems. In the time since Karp's paper, interest in NP-complete problems has exploded. We refer interested readers to Papadimitriou's book \cite{papa} on the topic.

The interest in NP-complete problems stems from the open question of whether there exists any polynomial-time algorithms to solve them. Since any $\mathcal{NP}$ problem can be reduced to any NP-complete problem, a polynomial-time algorithm for even one NP-complete problem would prove the existence of polynomial-time algorithms for all $\mathcal{NP}$ problems. However, to date no such algorithm has been found. This question is captured in the famous {\em P vs NP} problem, which asks whether $\mathcal{NP}$, and the set $\mathcal{P}$ of decision problems which are decidable in polynomial time, are equivalent.

Irrespective of whether a polynomial-time algorithm could be discovered for NP-complete problems, there is a second concern which, to date, has been largely ignored. Specifically, the question of how large the resultant instance is after a reduction is performed. Although it is, by definition, bounded above by a polynomial, the leading coefficient and the order of the polynomial may be large. Indeed, even if they are relatively low, the reduction from one problem to another may require several intermediate reductions, which compounds the input size of the final instance. In this context, the efficiency of performing the reductions is less important than the resultant size, as the total time taken to perform the reductions is merely the sum of the individual reduction times.

If $P \neq NP$, all NP-complete problems have exponential solving time in the worst case, and the merit of reducing any one of them to another problem is therefore lost if the input size grows too rapidly. Reductions usually only need to be performed for relatively large instances, since smaller instances can typically be solved by existing (exponential-time) algorithms. Even if $P = NP$, it is vital that a reduction does not inflate the input size of the instance too dramatically. Consider the following extremely optimistic prospect: a polynomial-time algorithm is discovered for an NP-complete problem that is guaranteed to conclude after $n^3$ iterations. Suppose you then wish to solve another NP-complete problem, but the reduction results in quadratic growth in input size. Even for a very modestly sized problem, say $n = 1000$, the resultant instance would take roughly $10^{18}$ iterations to solve, which is likely to be impractical. Hence, the order of the polynomial that bounds the input size of the new instance should be as small as possible; ideally, the polynomial should be linear, or at worst quasi-linear.

The above argument is the motivation for introducing the following definition.

\begin{definition}Consider two $\mathcal{NP}$ problems $P$ and $Q$. If a reduction exists from $P$ to $Q$ such that the input size of $Q$ is bounded above by a linear function of the input size of $P$, then we say that $Q$ lies in the {\em linear orbit} of $P$.\end{definition}

Obviously, if $Q$ is in the linear orbit of $P$, and in turn $R$ is in the linear orbit of $Q$, then $R$ is also in the linear orbit of $P$, so the property is transitive. However, it is not necessarily commutative. For example, boolean satisfiability with three literals per clause (3-SAT) is known to be in the linear orbit of Hamiltonian cycle problem (HCP), but HCP is not known to be in the linear orbit of 3-SAT, and indeed, it seems unlikely that it is. For completeness, we say that a problem is in its own linear orbit.

Then, consider a subset of $\mathcal{NP}$ called $\mathcal{NP}_L$, defined in such a way that any problem in $\mathcal{NP}$ is in the linear orbit of at least one problem in $\mathcal{NP}_L$, and $\mathcal{NP}_L$ is minimal. Certainly it seems reasonable that research efforts should be primarily focused on $\mathcal{NP}_L$ problems since these are the problems with the most potential scope for practical use. Indeed, if an efficient algorithm is developed for a problem with a large linear orbit, then all of those problems within its linear orbit can leverage off this algorithm as well, without needing to be concerned with such explosive growth as the example given earlier. Of course, a natural question to ask is whether $\mathcal{NP}_L$ is finite. Alternatively, if $\mathcal{NP}_L$ is not finite, what proportion of $\mathcal{NP}$ does it occupy?

In this manuscript, we focus on a more modest question, as a case study: if we consider solely the set, $\mathcal{K}_{21}$, of Karp's 21 NP-complete problems, how small a kernel subset $S \subset \mathcal{K}_{21}$ can we identify which possesses the property that all 21 problems lie in the linear orbit of at least one problem in $S$? The 21 NP-complete problems described by Karp are as follows:

\vspace*{0.45cm}\begin{minipage}[h]{0.45\linewidth}(1) Boolean satisfiability in conjunctive normal form (SAT)\\
(2) 0-1 Integer Programming\\
(3) Clique\\
(4) Set Packing\\
(5) Vertex Cover\\
(6) Set Covering\\
(7) Feedback Node Set\\
(8) Feedback Arc Set\\
(9) Directed Hamiltonian cycle problem (DHCP)\\
(10) Undirected Hamiltonian cycle problem (HCP)\end{minipage} \;\;\;
\begin{minipage}[h]{0.45\linewidth}(11) SAT with at most 3 literals per clause (3-SAT)\\
(12) Chromatic Number\\
(13) Clique Cover\\
(14) Exact Cover\\
(15) Hitting Set\\
(16) Steiner Tree\\
(17) 3-Dimensional Matching\\
(18) Knapsack\\
(19) Job Sequencing\\
(20) Partition\\
(21) Max Cut\end{minipage}

We provide reductions to demonstrate that $S$ can be reduced to cardinality six, specifically problems (2), (7), (10), (12), (13) and (19). We also discuss an ambiguity in the definition of input size that makes it unclear whether we should consider (12) and (13) as lying in each other's linear orbit. If so, then $S$ can be reduced to cardinality five. From an optimisation perspective, it is natural to assume that optimisation versions of problems in $S$ (for instance, see the discussion of complexity in Cook et al \cite{wcook}) may also be good surrogates for the optimisation versions of problems in their respective linear orbits.

\section*{Dealing with inequality constraints in 0-1 Integer Programming}

Throughout this manuscript, the majority of conversions will be to 0-1 Integer Programming. Using the definition given by Karp, only equality constraints are permitted. However, it will often be convenient to use inequality constraints. Of course, it is always possible to reduce an inequality-constrained integer program to an equality-constrained integer program through the use of slack and surplus variables. However, since we only permit binary variables, sometimes many slacks and surpluses will be needed. It is important to check very carefully how many slacks are required to ensure we do not exceed linear growth in the input size.

Consider the following example:
$$\sum_{i=1}^k x_i \leq m,$$
where $k,m$ are positive integers. Assuming that $x_i$ are binary variables, it is clear the LHS must be between 0 and $k$. When converting this constraint to an equality constraint, we must first ask ourself: what is the maximum difference between the LHS and RHS for which the inequality is still satisfied? It is clear that this situation occurs when the LHS is 0. Then we need to use as many slack variables $s_j$ as necessary to handle this situation. Define $n = \lceil \log_2(m+1) \rceil$. Then we can rewrite the above constraint as:
$$\sum_{i=1}^k x_i + \sum_{j=1}^n 2^{j-1} s_j = m.$$
It is easy to check that this constraint can be satisfied if and only if $x_i$ are chosen to satisfy the original inequality constraint. In the process of converting to an equality constraint, we introduced $n$ new variables, and $n$ new non-zero entries into the constraints coefficients matrix. The non-zero entries are $1, 2, 4, 8, \hdots, 2^{n-1}$ respectively. These can be encoded in binary using $1, 2, 3, 4, \hdots, n$ bits respectively. Hence, converting such an inequality constraint to an equality constraint increases the input size by $O(n^2) = O((\log_2(m))^2)$.

For problems where $m$ can grow with the size of the problem, care needs to be taken to ensure that this has not rendered the conversion super-linear. For example, suppose the input size of the original problem is $s$. If $m = O(2^s)$ then the above conversion produces constraints that require $O(s^2)$ memory to encode, and hence is not linearly-growing. Likewise, if $m = O(s)$ and there are $O(s)$ inequality constraints like the above, then after converting we require $O(s(\log_2(s))^2)$ memory to encode them all, and hence it also is not linearly-growing.

In the conversions that follow, we will consider these situations on a case-by-case basis to confirm that no such issues arise. Obviously, the same argument as above can be made when converting greater-than inequality constraints, using surplus variables, as well.

Note that if the maximum difference between the LHS and RHS in a valid solution is bounded above by a fixed constant, then the input size is increased by $O(1)$, and therefore does not prevent the conversion from being linearly-growing in any situation. In such a case we will say that the inequality constraint is {\em constant-bounded}.

\section*{Linearly-Growing Reductions}

Unless otherwise stated, all reductions in this paper are original. Those reductions which are not original all come from Karp's paper.

\subsection*{Satisfiability to 3-SAT (Karp)}

Satisfiability: Can a set of literals be assigned values of TRUE or FALSE so as to satisfy a set of clauses?\\
Input: $n$ clauses and $m$ literals. Each clause $C_i$ is of size $|C_i|$.\\
Input size: $\sum_i |C_i|$

3-SAT: Can a set of literals be assigned values of TRUE or FALSE so as to satisfy a set of clauses which all have cardinality 3?

Conversion: Produce a new instance of 3-SAT by constructing new clauses for each $C_i = \{\sigma_1, \sigma_2, \hdots, \sigma_{k}\}$, where $k = |C_i|$, as follows:

If $k \leq 3$: Simply repeat $C_i$.\\
If $k > 3$: Introduce new variables $y_{i1}, y_{i2}, \hdots, y_{i,k-3}$ and construct new clauses: $\hat{C}_{i1} = \{\sigma_1, \sigma_2, y_{i1}\},$ $\hat{C}_{i2} = \{\overline{y_{i1}}, \sigma_3, y_{i2}\},$ $\hat{C}_{i3} = \{\overline{y_{i2}}, \sigma_4, \overline{y_{i3}}\},$ $\hdots,$ $\hat{C}_{i,k-2} = \{\overline{y_{i,k-3}}, $ $\sigma_{k-1}, \sigma_k\}$.

Explanation: The intention is for the set new clauses to all be satisfiable if and only if the original clause $C_i$ is satisfiable. Consider the case where one of the literals in $C_i$, say $\sigma_n$, is assigned a value of TRUE. Then it can be seen that assigning $y_{i1} = y_{i2} = \hdots = y_{i,n-2} = $ TRUE, and $y_{i,n-1} = y_{in} = \hdots = y_{i,k-3} = $ FALSE satisfies all the clauses.

Next, consider the case where all of the literals in $C_i$ are assigned a value of FALSE. From $\hat{C}_{i1}$ it is clear that $y_{i1}$ must be assigned a value of TRUE. However, then from $\hat{C}_{i2}$, we see that $y_{i2}$ must also be assigned a value of true. Inductively, it follows that $y_{ij} = $ TRUE for $j = 1, \hdots, k-3$. However, this implies that clause $\hat{C}_{i,k-2}$ is not satisfied. Hence we conclude that the new clauses $\hat{C}_{ij}$ are all satisfiable if and only if the original clause $C_i$ was satisfiable.

Final Input Size: Consider each clause $C_i$ with $|C_i| > 3$. It is clear that in the new instance this has been replaced with $|C_i| - 2$ new clauses which each contain 3 literals. This implies that each clause $C_i$ (with input size $|C_i|$) has been replaced by new clauses with total input size $3|C_i| - 6$. Therefore an upper bound on the input size of the converted problem is $3 \sum_i |C_i|$. It is clear that the input size of the converted problem is a linear function of the original input size.

\subsection*{Clique to 0-1 Integer Programming}

Clique: Does a graph $G$ contain a clique (ie set of mutually adjacent vertices) of size $k$?\\
Input: Graph containing $N$ vertices and $e$ edges. Positive integer $k \leq N$.\\
Input size: $e + 1$.

0-1 Integer Programming: Is it possible to satisfy a set of linear equations in binary variables?

Conversion: Produce a new instance of 0-1 Integer Programming by introducing binary variables:

\begin{itemize}\item $x_{ij}$ for each edge $(i,j)$ in the graph,
\item $v_i$ for each vertex in the graph.\end{itemize}

Then the 0-1 Integer Program is described by the following constraints:

\begin{eqnarray}x_{ij} - v_i - v_j & \geq & -1, \qquad \forall \mbox{ edges } (i,j),\label{clique-1}\\
x_{ij} - v_i & \leq & 0, \qquad \forall \mbox{ edges } (i,j),\label{clique-2}\\
x_{ij} - v_j & \leq & 0, \qquad \forall \mbox{ edges } (i,j),\label{clique-3}\\
\sum_{i=1}^N v_i & = & k,\label{clique-4}\\
\sum_{(i,j)} x_{ij} & = & \left(\begin{array}{c}k\\2\end{array}\right).\label{clique-5}\end{eqnarray}

Explanation: The intention is for variables $v_i$ to be 1 if vertex $v$ is included in the $k$-clique, and 0 otherwise. Likewise, variables $x_{ij}$ will be 1 if edge $(i,j)$ is included in the $k$-clique, and 0 otherwise. Constraints (\ref{clique-1})--(\ref{clique-3}) ensure that $x_{ij} = 1$ if and only if both $v_i = 1$ and $v_j = 1$. Therefore the variables associated with edges between vertices in the clique are set to 1. Constraint (\ref{clique-4}) ensures that exactly $k$ vertices are in the clique, and constraint (\ref{clique-5}) ensures that there are sufficiently many edges between those vertices to constitute a clique.

Final Input Size: It can be checked that the constraints coefficients matrix will contain $N + 8e$ non-zero entries, and the RHS will contain $3e + 2$ entries. Note also that the RHS of constraint (\ref{clique-5}) will be less than $N^2$ and can therefore be encoded in at most $2\log(N)$ bits. All of the inequality constraints are constant-bounded. It is clear that the input size of the converted problem is a linear function of the original input size.

\subsection*{Set Packing to 0-1 Integer Programming}

Set Packing: Is it possible to select $l$ mutually disjoint sets from a family of sets?\\
Input: Family of sets $\{S_i\}$, with $s$ total sets and each set containing $|S_i|$ entries from a universe set $U$. Positive integer $l$.\\
Input size: $1 + \sum_{i=1}^s |S_i|$

0-1 Integer Programming: Is it possible to satisfy a set of linear equations in binary variables?

Conversion: Produce a new instance of 0-1 Integer Programming by introducing binary variables:

\begin{itemize}\item $x_i$ for $i = 1, \hdots, s$.\end{itemize}

Define $\sigma_{ij}$ to be 1 if set $S_i$ contains $j \in U$, and 0 otherwise. Then the 0-1 Integer Program is described by the following constraints:

\begin{eqnarray}\left(\sum_{i=1}^s \sigma_{ij} x_i \right) & \leq & 1, \qquad \forall j \in U,\label{setpacking-1}\\
\sum_{i=1}^s x_{i} & = & l.\label{setpacking-2}\end{eqnarray}

Explanation: The intention is for variables $x_i$ to be 1 if set $S_i$ is one of the $l$ mutually disjoint sets, and 0 otherwise. Constraint (\ref{setpacking-2}) ensures exactly $l$ sets are selected. Constraints (\ref{setpacking-1}) ensure that each entry of $U$ appears no more than once in the selected sets.

Final Input Size: It is easy to check that the number of non-zeros in the constraints coefficients matrix will be $\sum_{i=1}^s |S_i| + s$, and the RHS will contain $|U| + 1$ entries. Constraints (\ref{setpacking-1}) are constant-bounded inequalities. It is clear that the input size of the converted problem is a linear function of the original input size.

\subsection*{Node Cover to Set Covering (Karp)}

Node Cover: Is it possible to select no more than $l$ nodes in a graph $G$ such that every edge in $G$ is incident with at least one of the selected nodes?\\
Input: Graph containing $N$ vertices and $e$ edges. Positive integer $l \leq N$.
Input size: $e + 1$.

Set Covering: Is it possible to select no more than $k$ sets from a family of sets $F$, such that the union of the selected sets is equal to the union of all sets in $F$?

Conversion: Produce a new instance of Set Covering, in which each set contains elements taken from the set of edges in $G$ in the following manner. For each $j = 1, \hdots, N$, the set $S_j$ contains the edges incident with node $j$. Finally, assign $k$ (the constant in Set Covering) to be equal to $l$.

Explanation: Each set corresponds to a node in $G$. Since $k = l$, we can only select as many sets as we can select nodes. Then once $k$ sets are selected, a set covering is obtained if and only if every element in the sets is now covered. This corresponds to the situation where every edge in the graph is incident with at least one of the selected nodes.

Final Input Size: The number of entries over all of the sets will be $2e$, and the one constant input is precisely equal to $l$. It is clear that the input size of the converted problem is a linear function of the original input size.

\subsection*{Set Covering to 0-1 Integer Programming}

Set Covering: From a family of sets, is it possible to select no more than $k$ sets such that their union is equal to the union of all sets in the family?\\
Input: Family of sets $\{S_i\}$, with $s$ total sets and each set containing $|S_i|$ entries from a universe set $U$. Positive integer $k$.\\
Input size: $1 + \sum_{i=1}^s |S_i|$

0-1 Integer Programming: Is it possible to satisfy a set of linear equations in binary variables?

Conversion: Produce a new instance of 0-1 Integer Programming by introducing binary variables:

\begin{itemize}\item $x_i$ for $i = 1, \hdots, s$.\end{itemize}

Define $\sigma_{ij}$ to be 1 if set $S_i$ contains $j \in U$, and 0 otherwise. Then the 0-1 Integer Program is described by the following constraints:

\begin{eqnarray}\left(\sum_{i=1}^s \sigma_{ij} x_i \right) & \geq & 1, \qquad \forall i = 1, \hdots, s,\label{setcovering-1}\\
\sum_{i=1}^s x_i & = & k.\label{setcovering-2}\end{eqnarray}

Explanation: Although we only require at most $k$ sets, the problem generalises to choosing exactly $k$ sets. Then, the intention is for variables $x_i$ to be 1 if set $S_i$ is to be one of the $k$ sets chosen, and 0 otherwise. Constraint (\ref{setcovering-2}) ensures exactly $k$ sets are selected. Constraints (\ref{setcovering-1}) ensure that each entry of $U$ appears in at least one of the selected sets.

Final Input Size: It is easy to see that the number of non-zeros in the constraints coefficients matrix will be $\sum_{i=1}^s |S_i| + s$, and the RHS will contain $s + 1$ entries. Constraints (\ref{setcovering-1}) are not constant-bounded, so we consider them individually. For a particular choice of $i$, the difference between constraint the LHS and RHS could be as large as $|S_i|-1$. It is then clear that converting all $s$ constraints in (\ref{setcovering-1}) to equality constraints increases the size of the conversion by $O\left(\sum_{i=1}^s (\log_2(|S_i|))^2\right)$, which is certainly no bigger than the input size. It is clear that the input size of the converted problem is a linear function of the original input size.

\subsection*{Feedback Arc Set to Feedback Node Set}

Feedback Arc Set: Given a directed graph, is it possible to select a set of no more than $k$ arcs such that every (directed) cycle in the graph travels through at least one of the selected arcs?\\
Input: Graph containing $N$ vertices and $e$ (directed) arcs. Positive integer $k$.\\
Input size: $e + 1$.

Feedback Node Set: Given a directed graph, is it possible to select a set of no more than $k$ nodes such that every (directed) cycle in the graph travels through at least one of the selected nodes?

Conversion: We first convert the instance to an equivalent instance with nicer properties. Define the \lq\lq path graph" $P_n$ to be a graph containing $n$ vertices and arcs $(i,j)$ for every pair $i,j$ satisfying $|i - j| = 1$. Now consider the original graph, say $G$, and expand it in the following way. Suppose vertex $i$ has in-degree $c_i$ and out-degree $d_i$. Then replace vertex $i$ with $P_{c_i+d_i}$, ensuring that each arc incident on $i$ is now incident on a unique vertex of $P_{c_i+d_i}$. If there is a feedback arc set of size no more than $k$ in the original graph, there will be an equivalent one in this new graph, and vice versa. The new graph has bounded in-degree and out-degree of 3, and no more than $2e$ vertices. Hence there will be no more than $12e$ arcs in the new graph. Call the new graph $G'$.

Then the line graph of $G'$ constitutes an instance of Feedback Node Set with the identical choice of $k$.

Explanation: It is clear that a feedback node set in the line graph corresponds to a feedback arc set in $G'$. The reason the conversion to $G'$ is performed first is to obtain a sparse graph. This ensures the line graph is also sparse, and has size which is a linear function of the number of edges in $G'$, which in turn is linear in $e$. Hence we only need to show that $G'$ has a feedback arc set of size no more than $k$ if and only if $G$ does.

It is clear that we can obtain a feedback arc set of $G'$ from any feedback arc set of $G$ by simply selecting the corresponding arcs in $G'$, so the proof in one direction is trivial. Now consider the other direction. We will now view a feedback arc set as a set of arcs that may be removed, leaving a directed acyclic graph. Hence we can restrict our consideration to the cycles in $G'$. The only \lq\lq new" cycles in $G'$ (ie those that don't have a corresponding cycle in $G$) are those that are created by being allowed to visit one or more path graphs multiple times. These effectively correspond to a union of cycles in $G$. Hence if there is no feedback arc set of size $k$ that removes all the cycles in $G$, there is definitely none of size $k$ in $G'$ either.

Final Input Size: As argued above, the line graph of $G'$ will be sparse (ie have in-degree and out-degree bounded by a constant) and will have a vertex set of cardinality a linear function of $e$. It is clear that the input size of the converted problem is a linear function of the original input size.

\subsection*{Directed HCP to Undirected HCP (Karp)}

Directed HCP: Does a given directed graph contain a simple cycle that visits every vertex?\\
Input: Graph $G$ containing $N$ vertices and $e$ (directed) arcs.\\
Input size: $e$.

Undirected HCP: Does a given undirected graph contain a simple cycle that visits every vertex?

Conversion: Define subgraphs $S_i$ for $i = 1, \hdots, N$ to be 3-vertex subgraphs each containing edges $(1,2)$ and $(2,3)$. We then construct a new instance of Undirected HCP by replacing each vertex $i$ in $G$ with $S_i$. Then for each directed arc $(i,j)$, the new instance contains an edge going from vertex 3 of $S_i$ to vertex 1 of $S_j$.

Explanation: The second vertex in each $S_i$ is a degree 2 vertex, and so it is clear that any time vertex $1$ of any $S_i$ is reached, vertices 2 and then 3 of the same $S_i$ must immediately follow. Then it is only possible to exit each $S_i$ via an edge incident on the third vertex, which corresponds to an arc that departs vertex $i$ in $G$. Likewise, each time a $S_i$ is exited, another $S_j$ is entered via vertex 1, which corresponds to an arc that enters vertex $j$ in $G$. Then it is clear that Hamiltonian cycles in $G$ have a 1-1 relationship with the Hamiltonian cycles in the converted instance.

Final Input Size: Each directed arc in $G$ now has a corresponding edge in the converted instance. In addition, there are two extra edges for each $S_i$. It is clear that the input size of the converted problem is a linear function of the original input size.

\subsection*{3-SAT to 0-1 Integer Programming (Karp)}

3-SAT: Can a set of literals be assigned values of TRUE or FALSE so as to satisfy a set of clauses $C_i$, where $|C_i| = 3$ for all $i$?\\
Input: $n$ clauses and $m$ literals. Each clause $C_i$ is of size $3$.\\
Input size: $3n$

0-1 Integer Programming: Is it possible to satisfy a set of linear equations in binary variables?

Conversion: Produce a new instance of 0-1 Integer Programming by introducing binary variables:

\begin{itemize}\item $x_j$ for $j = 1, \hdots, m$.\end{itemize}

Suppose the entries in clause $C_i$ are $(c_{i1}, c_{i2}, c_{i3})$. Define $d_i$ to be the number of complemented variables in $C_i$, and define $w_{ij}$ as follows:

$$w_{ij} = \left\{\begin{array}{rcl}1 & \mbox{ if } & c_{ik} = x_j \mbox{ for any } k = 1, 2, 3\\
-1 & \mbox{ if } & c_{ik} = \bar{x}_j \mbox{ for any } k = 1, 2, 3\\
0 & & \mbox{otherwise}\end{array}\right.$$

Then the 0-1 Integer Program is described by the following constraints:

\begin{eqnarray}\left(\sum_{j=1}^m w_{ij}x_j\right) & \geq & 1 - d_i, \qquad \forall i = 1, \hdots, n.\label{3sat-1}\end{eqnarray}

Explanation: The intention is for variables $x_j$ to be 1 if literal $i$ is to be assigned TRUE, and 0 if literal $i$ is to be assigned FALSE. Then for each clause $i$ we want at least one of the literals to have the desired value. If the literal $j$ is not complemented, we include the variable $x_j$. If it is complemented, then we include $(1 - x_j)$. Rearranging this, we see that for each clause, we must satisfy constraint (\ref{3sat-1}). Likewise, if constraint (\ref{3sat-1}) is satisfied, then there is a valid assignment of literals that satisfies all of the clauses.

Final Input Size: The number of non-zeros in the constraints coefficients matrix is $3n$, and there are $n$ RHS entries. Constraints (\ref{3sat-1}) are all constant-bounded. It is clear that the input size of the converted problem is a linear function of the original input size.

\subsection*{Exact Cover to 0-1 Integer Programming}

Exact Cover: Given a family of sets, is it possible to select a subfamily of mutually disjoint sets whose union is equal to the union of all sets in the family?\\
Input: Family of sets $\{S_i\}$, with $s$ total sets and each set containing $|S_i|$ entries from a universe set $U$.\\
Input size: $\sum_{i=1}^s |S_i|$.

0-1 Integer Programming: Is it possible to satisfy a set of linear equations in binary variables?

Conversion: Produce a new instance of 0-1 Integer Programming by introducing binary variables:

\begin{itemize}\item $x_i$ for $i = 1, \hdots, s$.\end{itemize}

Define $\sigma_{ij}$ to be 1 if set $S_i$ contains $j \in U$, and 0 otherwise. Then the 0-1 Integer Program is described by the following constraints:

\begin{eqnarray}\sum_{i=1}^s \sigma_{ij} x_i & = & 1, \qquad \forall j \in U.\label{exactcover-1}\end{eqnarray}

Explanation: The intention is for variables $x_i$ to be 1 if set $S_i$ is to be included in the exact cover, and 0 otherwise. Then constraints (\ref{exactcover-1}) ensure that each entry in $U$ appears precisely once in the selected sets.

Final Input Size: The number of non-zeros in the constraints coefficients matrix is precisely equal to $\sum_{i=1}^s |S_i|$, and there are $|U|$ RHS entries. It is clear that the input size of the converted problem is a linear function of the original input size.

\subsection*{Hitting Set to 0-1 Integer Programming}

Hitting Set: Given a family of sets, is it possible to construct a new set $W$ such that the intersection between the $W$ and any of the given sets has cardinality 1?\\
Input: A family of sets $\{S_i\}$, with $s$ total sets and each set containing $|S_i|$ entries from a universe set $U$.\\
Input size: $\sum_{i=1}^s |S_i|$.

0-1 Integer Programming: Is it possible to satisfy a set of linear equations in binary variables?

Conversion: Produce a new instance of 0-1 Integer Programming by introducing binary variables:
\begin{itemize}\item $x_j$ for $j = 1, \hdots, |U|$.\end{itemize}
Define $\sigma_{ij}$ to be 1 if set $S_i$ contains $j \in U$, and 0 otherwise. Then the 0-1 Integer Program is described by the following constraints:
\begin{eqnarray}\sum_{j=1}^{|U|} \sigma_{ij} x_j & = & 1, \qquad \forall i = 1, \hdots, s.\label{hittingset-1}\end{eqnarray}
Explanation: The intention is for variables $x_i$ to be 1 if they are to be included in $W$, and 0 otherwise. Then constraints (\ref{hittingset-1}) ensure that each set $S_i$ contains precisely one of the selected variables.

Final Input Size: The number of non-zeros in the constraints coefficients matrix is precisely equal to $\sum_{i=1}^s |S_i|$, and there are $s$ RHS entries. It is clear that the input size of the converted problem is a linear function of the original input size.

\subsection*{Steiner Tree to 0-1 Integer Programming}
Steiner Tree: Given a graph $G$, weights for each edge, and a set of vertices $R$ in the graph $G$, is it possible to find a subtree in $G$ containing all vertices in $R$ such that the total weight of the tree no more than a given value $k$?\\
Input: Graph containing $N$ vertices and $e$ edges. Set of vertices $R$. Set of weights $W$ made up of $e$ values of the form $w_{ij}$. Positive integer $k$.\\
Input size: $2e + |R| + 1$.

0-1 Integer Programming: Is it possible to satisfy a set of linear equations in binary variables?

Conversion: Then, produce a new instance of 0-1 Integer Programming by introducing binary variables:

\begin{itemize}\item $x_{ij}$ for each edge $(i,j)$,
\item $y_j$ for $j = 1, \hdots, N$,
\item $z_j$ for $j = 1, \hdots, N$.\end{itemize}

Define $\mathcal{A}(j)$ to be the set of vertices adjacent to $j$, that is $i \in \mathcal{A}(j)$ if and only if edge $(i,j)$ exists in $G$. Then the 0-1 Integer Program is described by the following constraints:

\begin{eqnarray}\sum_{j=1}^N z_j & = & 1,\label{steinertree-1}\\
y_j + z_j & = & 1, \qquad \forall j \in R,\label{steinertree-2}\\
y_j + z_j & \leq & 1, \qquad \forall j \not\in R,\label{steinertree-3}\\
y_j - \sum_{i \in \mathcal{A}(j)} x_{ij} & = & 0, \qquad \forall j = 1, \hdots, N,\label{steinertree-4}\\
x_{ij} - y_i - z_i & \leq & 0, \qquad \forall \mbox{ edges } (i,j),\label{steinertree-5}\\
x_{ij} + z_j & \leq & 1, \qquad \forall \mbox{ edges } (i,j),\label{steinertree-6}\\
\left(\sum_{\mbox{edges }(i,j)} w_{ij} x_{ij}\right) & \leq & k.\label{steinertree-7}\end{eqnarray}

Explanation: The intention is for variables $x_{ij}$ to be 1 if edge $(i,j)$ is to be used in the subtree, and 0 otherwise. The $x_{ij}$ variables are oriented in the sense that edges must emanate from lower in the tree, so for example, if vertex $i$ is the root vertex and edge $(i,j)$ appears in the tree, then $x_{ij} = 1$ but $x_{ji} = 0$. Variable $z_j$ is designed to be 1 if vertex $j$ is the root of the subtree, and 0 otherwise. Similarly, variable $y_j$ is designed to be 1 if vertex $j$ is any member of the subtree other than the root vertex, and 0 otherwise.

Constraint (\ref{steinertree-1}) ensures there is only one root. Constraints (\ref{steinertree-2}) ensure that every vertex in $R$ is either the root of the subtree, or another member of the subtree. Constraints (\ref{steinertree-2})--(\ref{steinertree-3}) ensure that no vertex is viewed as being both the root of the subtree and also another member of the subtree. Constraints (\ref{steinertree-4}) ensure that vertices are only seen as (non-root) members of the subtree if a single edge enters it (as is the definition of a tree). Constraints (\ref{steinertree-5}) ensure that any edge $(i,j)$ may only be used if vertex $i$ is a member of the subtree. Constraints (\ref{steinertree-6}) ensure that any edges incident on the root vertex and contained in the subtree must only emanate from the root vertex, rather than go to the root vertex.

Constraints (\ref{steinertree-1})--(\ref{steinertree-6}) combine to ensure that the set of $x_{ij}$ variables correspond to a valid, connected tree that contains $R$. Then, finally, constraint (\ref{steinertree-7}) ensures that the weight of the subtree does not exceed $k$.

Final Input Size: It can be checked that the number of non-zeros in the 0-1 Integer Program is $4N + 7e$. Note that $e > N$ in any meaningful example. Inequality constraints (\ref{steinertree-3}), (\ref{steinertree-5}) and (\ref{steinertree-6}) are all constant-bounded, but constraint (\ref{steinertree-7}) is not. The maximum difference between the LHS and RHS of constraint (\ref{steinertree-7}) is $k$, which theoretically can grow infinitely large. However, if $k$ is larger than the sum of all weights in the graph, then constraint (\ref{steinertree-7}) is satisfied automatically and can be ignored. Hence, we assume that $k$ is not larger than the sum of weights, and so converting (\ref{steinertree-7}) to an equality constraint will certainly increase the problem size by less than the original input size. It is clear that the input size of the converted problem is a linear function of the original input size.

\subsection*{3-Dimensional Matching to 0-1 Integer Programming}

3-Dimensional Matching: Given $U \subseteq T \times T \times T$, is it possible to find $W \subseteq U$ such that $W$ contains $|T|$ entries from $U$, and no two entries of $W$ agree in any coordinate?\\
Input: A family $U$, containing $|U|$ sets (the $i$-th set being called $U_i$) which each contain three entries. The finite size $|T|$.\\
Input size: $3|U| + 1$

0-1 Integer Programming: Is it possible to satisfy a set of linear equations in binary variables?

Conversion: Produce a new instance of 0-1 Integer Programming by introducing binary variables:

\begin{itemize}\item $x_i$ for $i = 1, \hdots, |U|$.\end{itemize}

Define $\sigma_{ijk}$ to be 1 if set $U_i$ contains entry $j \in T$ in coordinate $k$, and 0 otherwise. Then the 0-1 Integer Program is described by the following constraints:

\begin{eqnarray}\sum_{i=1}^{|U|} \sigma_{ijk} x_i & = & 1, \qquad \forall j = 1, \hdots, |T|, \;\;\; k = 1, 2, 3.\label{3dmatching-1}\end{eqnarray}

Explanation: The intention is for variables $x_i$ to be 1 if set $U_i$ is included in $W$, and 0 otherwise. Since there should be $|T|$ sets included in $W$, and no two entries of $W$ are to agree in any coordinate, it is clear that the sets of $W$ will cover every single entry in $T$ for all three coordinates precisely once. Conversely, if every entry in each coordinate appears precisely once, then it must be the case that $|W| = |T|$, as desired. To that end, constraints (\ref{3dmatching-1}) request that each entry appears precisely once, which will only be possible if a 3-dimensional matching can be found.

Final Input Size: The number of non-zeros in the 0-1 Integer Program is precisely equal to $3|U|$, and there are $3|T|$ RHS entries. For any meaningful instance of 3-Dimensional Matching, $|U| > |T|$. It it clear that the input size of the converted problem is a linear function of the original input size.

\subsection*{Knapsack to 0-1 Integer Programming}

Knapsack: Given a set of integers and a target value $b$, is it possible to choose some integers from the set so that their sum is equal to $b$?\\
Input: A set $A$ containing $r$ integers, and a target value $b$.\\
Input size: $r + 1$

0-1 Integer Programming: Is it possible to satisfy a set of linear equations in binary variables?

Conversion: Produce a new instance of 0-1 Integer Programming by introducing binary variables:

\begin{itemize}\item $x_i$ for $i = 1, \hdots, r$.\end{itemize}

Denote by $a_i$ the $i$-th entry of $A$. Then the 0-1 Integer Program is described by the following constraints:

\begin{eqnarray}\sum_{i=1}^{r} a_i x_i & = & b.\label{knapsack-1}\end{eqnarray}

Explanation: The intention is for variables $x_i$ to be 1 if integer $a_i$ is chosen in the sum, and 0 otherwise. Then, it is clear by definition that the sole constraint (\ref{knapsack-1}) describes the Knapsack problem perfectly.

Final Input Size: There are exactly $r$ non-zero entries in the 0-1 Integer Program, and a single RHS entry. Then it is clear that the input size of the converted problem is precisely equal to the input size of the original problem.

\subsection*{Partition to Knapsack}

Partition: Given a set of integers, is it possible to choose some integers from the set so that their sum is equal to the sum of the integers not selected?\\
Input: A set $C$ containing $s$ integers.\\
Input size: $s$

Knapsack: Given a set of integers and a target value $b$, is it possible to choose some integers from the set so that their sum is equal to $b$?\\

Conversion: Denote by $c_i$ the $i$-th entry of $C$. Then the new instance of Knapsack is produced by simply $C$ as the new set of integers, and setting $b = \frac{1}{2}\sum_{i=1}^s c_i$.

Explanation: Clearly, if integer can be chosen from $C$ such that the total is half of the sum of all integers in $C$, then the remaining integers will also sum up to the same value, satisfying the Partition condition.

Final Input Size: Since $C$ is provided for both Partition and Knapsack, the only additional input is $b = \frac{1}{2}\sum_{i=1}^s c_i$. Since $b$ is a number less than the sum of all entries of $C$, it can certainly be encoded to be no larger than the input of $C$ itself. It is clear that the input size of the converted problem is a linear function of the original input size.

%
%
%
%
%
%
%
%

\subsection*{Max Cut to 0-1 Integer Programming}

Max Cut: Given a graph $G$ with weights on each edge, and a positive integer $W$, is it possible to select a set of vertices $S$ such that the sum of weights on edges with exactly one vertex in $S$ is at least as big as $W$?\\
Input: Graph $G$ containing $N$ vertices and $e$ edges. Weights $w_{ij}$ for all $e$ edges. Positive integer $W$.\\
Input Size: $2e + 1$

0-1 Integer Programming: Is it possible to satisfy a set of linear equations in binary variables?

Conversion: Define $TW = \sum_{(i,j)} w_{ij}$. Then produce a new instance of 0-1 Integer Programming by introducing binary variables;

\begin{itemize}\item $x_i$ for $i = 1, \hdots, N$,
\item $y_{ij}$ for each edge $(i,j)$.\end{itemize}

Then the 0-1 Integer Program is described by the following constraints:


\begin{eqnarray}y_{ij} - x_i + x_j & \leq & 1, \qquad \forall \mbox{ edges } (i,j),\label{maxcut-1}\\
y_{ij} - x_j + x_i & \leq & 1, \qquad \forall \mbox{ edges } (i,j),\label{maxcut-2}\\
y_{ij} + x_i + x_j & \geq & 1, \qquad \forall \mbox{ edges } (i,j),\label{maxcut-3}\\
y_{ij} - x_i - x_j & \geq & -1, \qquad \forall \mbox{ edges } (i,j),\label{maxcut-4}\\
\left(\sum_{(i,j)} w_{ij} y_{ij}\right) & \leq & TW - W.\label{maxcut-5}\end{eqnarray}

Explanation: The intention is for variables $x_i$ to be 0 if vertex $i$ is to be selected in $S$, and 1 otherwise. Also, $y_{ij}$ is to be set to 0 if edge $(i,j)$ is added in the sum, and 1 otherwise. The variables are chosen this way so the problem can be reformulated with a less-than inequality rather than greater-than. To that end, constraints (\ref{maxcut-1})--(\ref{maxcut-4}) are designed in such a way that they can only all be satisfied if the following condition is true: $y_{ij} = 1$ if and only if $x_i + x_j \neq 1$. Finally, constraint (\ref{maxcut-5}) ensures that the weights of all edges that {\em don't} have exactly one vertex in $S$ is no bigger than $TW - W$, which is equivalent to the desired condition on $W$.

Final Input Size: It can be checked that the number of non-zeros in the 0-1 Integer Program is $13e$, and there are $4e + 1$ RHS entries. Constraints (\ref{maxcut-1})--(\ref{maxcut-4}) are all constant-bounded. Constraint (\ref{maxcut-5}) is not constant-bounded, and the RHS is a potentially large number not used as input in the original problem, so we consider this constraint individually. Consider first encoding the number $TW - W$. This number will be no larger than the sum of all weights, and so it can be encoded in fewer bits than it takes to encode the weights in the original input problem. Then, it is clear that the maximum difference between the LHS and RHS of (\ref{maxcut-5}) is no larger than $TW$, so the increase in problem size after converting (\ref{maxcut-5}) to an equality constraint is smaller than the original input size as well. It is clear that the input size of the converted program is a linear function of the original input size.

\section*{Ambiguity in Input Size}

In the previous section we have determined the input size of each problem by considering the amount of data required to store the information that describes the problem instance. However, it is not always clear how this should be computed. For example, is it reasonable to assume a problem is stored in a format that requires preprocessing to be used by an algorithm? Alternatively, is it reasonable to think of input size as the amount of information required to be stored in memory by a standard algorithm for that problem? We highlight this ambiguity with the following example.

\subsection*{Chromatic Number to Clique Cover (Karp)}

Chromatic Number: Is it possible to assign a colour to each vertex of a graph such that no two adjacent vertices have the same colour, and the total number of colours is no bigger than $k$?\\
Input: A graph $G$ containing $N$ vertices and $e$ edges. A positive integer $k \leq N$.\\
Input size: $e + 1$

Clique Cover: Is it possible to select no more than $l$ cliques in a graph such that none of the cliques overlap and the union covers all vertices?

In the following, we provide a conversion for which the input size is ambiguously defined.

Conversion: Produce a new instance of Clique Cover by defining:

$G' = $ the complement of $G$\\
$l = k$.

Explanation: Suppose there is a clique cover of $G'$ containing cliques $C_i$ for $i = 1, \hdots, j$ and $j \leq l$. Then consider the vertices in clique $C_i$. One can safely colour these vertices the same colour in $G$, as by definition of $G'$ none will be adjacent to any other in $G$. It is clear that the two problems are then equivalent.

Potential Issue: It is unclear how to define the input size of the resultant problem. Technically, if $G$ is a sparse graph, then $G'$ is a dense graph. In this case, the input size of the new Clique Cover instance would be $O(Ne)$. However, in this case it is possible to define $G'$ by storing a sparse amount of information (ie by storing $G$) and then taking the complement once $G$ has been read in. Should the size of the problem be described in terms of the most efficient method of storage, or the number of elements in the problem? If we imagine that we had a black box solver for Clique Cover that required us to input a graph directly, we would not be able to submit a different graph and demand that the solver first complement the graph. Alternatively, if we consider the input size to be the number of elements that will be manipulated by any potential solver, then we would be required to consider the dense graph explicitly.

If we consider the size of the problem being purely the number of bits it takes to encode the problem in its most compressed form, the above constitutes a linearly-growing reduction. However, if we do not permit special processing of the instance, then it does not.

If $G$ is a dense graph, then the above conversion is linearly-growing regardless of how we store the data.

\section*{Discussion and Future Work}

In this paper, we introduced the notion of a {\em linear orbit} of an problem NP-complete $P$. Namely, the orbit consists of the set of problems in $\mathcal{NP}$ that can be converted to $P$ by a conversion that results in linear growth in the input size. In particular, we showed that there exists a kernel subset $S$ of the set $\mathcal{K}_{21}$ of the classical NP-complete problems stated in Karp's seminal 1972 paper. Every one of the 21 problems belongs to the linear orbit of at least one of the six problems in $S$.

These results suggest that efficient algorithms to solve problems in $S$ may offer opportunities to more efficiently solve the problems in their linear orbits. It is hoped that this can be extended from problems in their decision framework to the more practical optimisation frameworks.

It is worth contemplating how much smaller $S$ would be if we were to permit reductions with larger growth. For example, if we permit reductions that result in quasilinear growth, it is possible to reduce HCP to SAT \cite{johnson}. If we further permit $O(n\sqrt{n})$ growth, the reduction of Chromatic Number to Clique Cover above is permitted, and it is also possible to show that Clique Cover could be reduced to 0-1 Integer Programming. If quadratic growth is permitted, Job Sequencing can be reduced to 0-1 Integer Programming. Clearly this is a topic ripe for future research.

\medskip
Received September 2015;  1$^{st}$ revision August 2016; final revision January  2018.
\medskip

\end{document}